\documentclass[twoside,10pt]{amsart}
\title{Weak covering properties and infinite games}
\author{L. Babinkostova, B.A. Pansera and M. Scheepers}
\date{}

\usepackage{amssymb}

\newtheorem{theorem}{Theorem}
\newtheorem{definition}[theorem]{Definition}

\newtheorem{question}[theorem]{Question}

\newtheorem{lemma}[theorem]{Lemma}
\newtheorem{corollary}[theorem]{Corollary}
\newtheorem{proposition}[theorem]{Proposition}

\def\Proof{\textbf{Proof.  }}
\newcommand{\open}{{\mathcal O}}
\newcommand{\dense}{{\mathcal D}}
\newcommand{\almost}{\overline{\open}}
\newcommand{\sone}{{\sf S}_1}
\newcommand{\sfin}{{\sf S}_{fin}}
\newcommand{\gone}{{\sf G}_1}
\newcommand{\gfin}{{\sf G}_{fin}}

\newcommand{\htheta}{{\sf H}_{\theta}}
\newcommand{\modelM}{{\mathcal M}}

\newcommand{\naturals}{{\mathbb N}}
\newcommand{\reals}{{\mathbb R}}

\newcommand{\cohen}{{\mathbb C}}
\newcommand{\poset}{{\mathbb P}}
\newcommand{\forces}{\mathrel{\|}\joinrel\mathrel{-}}

\address{Department of Mathematics\\ Boise State University\\ Boise, Idaho 83725}
\email{liljanababinkostova@boisestate.edu,\\
       mscheepe@boisestate.edu}
\address{Dipartimento di Matematica, Universita di Messina, Via F. Stagno d'Alcontres N.31,
          8166 Messina (Italy). (Pansera) }
\email{bpansera@unime.it}

\begin{document}

\maketitle

\begin{abstract} We investigate game-theoretic properties of selection principles related to weaker forms of the Menger and Rothberger properties. For appropriate spaces some of these selection principles are characterized in terms of a corresponding game. We use generic extensions by Cohen reals to illustrate the necessity of some of the hypotheses in our theorems. 
\end{abstract}

\section{Introduction} 

A \emph{Lindel\"of space} is a topological space for which each open cover contains a countable subset covers the space. Some unresolved or undecidable problems about Lindel\"of spaces have definite answers for spaces having a stronger version of the Lindel\"of property. For example: Cardinality questions that are undecidable for classes of Lindel\"of spaces have been resolved for the corresponding classes of \emph{compact} spaces. While it is not known if regular ${\sf T}_3$ Lindel\"of spaces are {\sf D}-spaces\footnote{As we do not need the concept of a {\sf D}-space elsewhere in our paper, it is left undefined.} it is known that each ${\sf T}_1$ \emph{Menger} space (a selective version of the Lindel\"of property and defined later in our paper) is a {\sf D}-space. Whereas the Lindel\"of property is not preserved by the product construction, compactness is.

In some cases where the Lindel\"of property (or one of its strengthenings) is not preserved by a topological construction, it has been found that a weaker version of the Lindel\"of property is preserved. Also, it has been found that some theorems true for Lindel\"of spaces are in fact true for a wider class of spaces that have a weak version of the Lindel\"of property. We consider two such weakenings of the Lindel\"of property that have emerged in the literature. 

To define these, let $\open$ denote the collection of all open covers of a space, let $\dense$ denote the collection of families of open sets with union dense in the space, and let $\almost$ denote the collection of families $\mathcal{U}$ of open subsets of the space for which $\{\overline{U}:U\in\mathcal{U}\}$ covers the space\footnote{The symbol $\overline{V}$ denotes the closure of the set $V$. The notation $\overline{{\mathcal O}}$ is due to Boaz Tsaban.}. A space is said to be \emph{weakly Lindel\"of} if each open cover contains a countable subset with union dense in the space. Thus, weakly Lindel\"of means that each element of $\mathcal{O}$ has a countable subset which is a member of $\mathcal{D}$. The notion of weakly Lindel\"of appears to have been introduced in the 1959 paper \cite{F} of Frolik. The name ``weakly Lindel\"of" however seems to have been coined by Hager and Mrowka (see the introductory paragraph of \cite{CHN}). A thorough introduction to weakly Lindel\"of spaces can be found in \cite{CHN}. A topological space $X$ is said to be \emph{almost Lindel\"of} if there is for each open cover $\mathcal{U}$ of $X$ a countable subset $\mathcal{V}$ such that $\{\overline{V}:V\in\mathcal{V}\}$ is a cover of $X$. The notion of an almost Lindel\"of space was introduced in the 1984 paper \cite{WD} by Willard and Dissanyake. The following implications hold among these three properties:

\begin{center}
Lindel\"of $\Rightarrow$ almost Lindel\"of $\Rightarrow$ weakly Lindel\"of.
\end{center}

Various selective versions of the Lindel\"of property, for example the Menger property or the Rothberger property (defined below), have characterizations in terms of infinite games. These game characterizations are powerful tools to derive other mathematical properties of these classes of spaces. We investigate the possibility of characterizing certain selective versions of the almost Lindel\"of and the weakly Lindel\"of properties by infinite games. To indicate to what extent some of the hypotheses of some of our results are necessary we also explore the preservation of selective versions of these properties under generic extensions of the universe. Behavior of these properties under topological constructions  will be addressed in the paper \cite{BPSProducts}.

\section{Definitions}

Recall the general framework for describing selection hypotheses (as in \cite{COC2} and \cite{COC1}):
Let $\bf N$ denote the set of positive integers. Let $\mathcal A$ and $\mathcal B$ be collections of subsets of an infinite set. Then
$\sone({\mathcal A}, {\mathcal B})$ denotes the following hypothesis:
\begin{quote}
For each sequence $(A_n : n \in \naturals)$ of elements of $\mathcal A$ there is a sequence $(B_n : n\in {\mathbb N})$ such that, for each $n$, $B_n\in A_n$ and $\{B_n : n \in \naturals\}$ is an element of $\mathcal B$.
\end{quote}
Then $\sone(\open,\open)$ denotes the \emph{Rothberger} property.

 The symbol $\sfin({\mathcal A}, {\mathcal B})$ similarly denotes the hypothesis
\begin{quote}
For each sequence $(A_n : n \in \naturals)$ of elements of $\mathcal A$ there is a sequence $(B_n : n\in \naturals)$ such that, for each $n$, $B_n\subseteq A_n$ is finite, and $\bigcup\{B_n : n \in \naturals\}$ is an element of $\mathcal B$.
\end{quote}
$\sfin(\open,\open)$ denotes the \emph{Menger} property.

Our conventions for the rest of the paper are: By ``space" we mean a topological space. Unless stronger separation axioms are
indicated specifically, we assume all spaces to be infinite and ${\sf T}_1$. Undefined notation and terminology will be as in \cite{E}.

\section{Preserving Lindel\"of like properties in Cohen real generic extensions}

It is well known that 
\begin{theorem}\label{nonLpreserve} Let $\poset$ be a proper partial order and let $X$ be a space.
\begin{enumerate}
  \item{If $X$ is non-Lindel\"of, then ${\mathbf 1}_{\poset}\forces``\check{X} \mbox{ is not Lindel\"of}"$.}
  \item{If $X$ is not almost Lindel\"of, then ${\mathbf 1}_{\poset}\forces``\check{X} \mbox{ is not almost Lindel\"of}"$.}
  \item{If $X$ is not weakly Lindel\"of, then ${\mathbf 1}_{\poset}\forces``\check{X} \mbox{ is not weakly Lindel\"of}"$.}  
\end{enumerate}
\end{theorem}
\Proof These statements follow directly from Proposition 4.1 of \cite{Jech}: Consider a ground model that  satisfies {\sf ZFC}. If $\poset$ is a proper partially ordered set in this ground model and $S$ is a ground model set, then for a countable set $Y$ that is a subset of $S$ but a member of the generic extension, there is a countable set $Z$ that is a subset of $S$ and is a member of the  ground model, and contains $Y$.  $\Box$

The Lindel\"of property is preserved by adding Cohen reals or random reals (see for example \cite{ST}), but need not be preserved by countably closed forcing \cite{Tall} or even certain ccc forcing extensions \cite{Isaac}. More information about these issues can be found in \cite{Kada}, \cite{ST} and \cite{Tall}.
We now establish some corresponding results for the almost Lindel\"of and the weakly Lindel\"of properties.  

For $\kappa$ a fixed uncountable cardinal $\cohen(\kappa)$ denotes the partially ordered set for adding $\kappa$ Cohen reals generically: The underlying set of elements of $\cohen(\kappa)$ is 
\[
  {\sf Fn}(\kappa\times\omega,\omega)=\{p\subset (\kappa\times\omega)\times \omega: p \mbox{ a finite function with domain a subset of }\kappa\times\omega\}
\]
For $p$ and $q$ in ${\sf Fn}(\kappa\times\omega,\omega)$ we write $p<q$ if $q\subset p$. It is well known that forcing with $\cohen(\kappa)$ preserves cardinals and cofinalities, and that in the resulting model $2^{\aleph_0}\ge\kappa$. Since $\cohen(\kappa)$ is a proper partially ordered set, Theorem \ref{nonLpreserve} applies to it.

\begin{definition}  For a positive integer $n$, an  n-dowment  is  a  family  $\mathcal{L}_n$  of  finite  antichains  of $\cohen(\kappa)$  such 
that: 
\begin{enumerate}
\item{For  each  maximal  antichain  $\mathbb{A}  \subset  \cohen(\kappa)$  there  is  an  $L  \in  \mathcal{L}_n$  such  that  $L  \subseteq \mathbb{A}$}
\item{For  each  $p \in  \cohen(\kappa)$  such  that  $\vert dom(p)\vert  <  n$  and  for  every  collection  $L_1,  L_2, \cdots L_n\in {\mathcal L}_n$ there  are  $q_1  \in  L_1,  \cdots  , q_n  \in  L_n$,  and  there  is  $r  \in  \cohen(\kappa)$,  such  that $r  \leq  p$  and  $r  \leq  q_i$  for  each  $i  \leq  n$.} 
\end{enumerate}
\end{definition}

\begin{center}{\bf Preserving weakly Lindel\"of}\end{center}

\begin{theorem}\label{cohenrealpreservation} Let $\kappa>\aleph_0$ be a cardinal. Let $X$ be a topological space.
  {If $X$ is weakly Lindel\"of, then in ${\sf V}^{\cohen(\kappa)}$ $X$ is weakly Lindel\"of.} 
\end{theorem}
$\Proof$
   For each $n$ let ${\mathcal L}_n$ be an $n$-dowment for $\cohen(\kappa)$.
   Let $\tau$ denote the (ground model) topology of $X$. Let $\dot{\mathcal{U}}$ be a $\cohen(\kappa)$-name such that
\[
  {\mathbf {1}_{\cohen(\kappa)}}\forces``\dot{\mathcal{U}}\subseteq \check{\tau} \mbox{ is an open cover of }\check{X}".
\]
  For each $x$ choose a maximal antichain ${\mathbb A}_x\subseteq \cohen(\kappa)$ and ground model open sets $V_{p,x}$, $p\in{\mathbb A}_x$, such that $x\in  V_{p,x}$ and $p \forces ``\check{V}_{p,x}\in \dot{U}"$.
For each $n$ choose ${\mathbb A}_{n,x}\in{\mathcal{L}}_n$ and define $V_{n,x} = \cap\{V_{p,x}: p\in {\mathbb A}_{n,x}\}$. Then $V_{n,x}$ is open, and $\mathcal{U}_n = \{V_{n,x}:x\in X\}$ is an open cover of $X$ in the ground model.

Since $X$ is weakly Lindel\"of, choose a countable set $\mathcal{V}_n\subseteq \mathcal{U}_n$ such that $\bigcup{\mathcal V}_n$ is dense in $X$. 

{\flushleft{\bf Claim: }} ${\mathbf 1}_{\cohen(\kappa)}\forces ``(\forall V\in\check{\tau})(\exists n)(\exists U\in\check{\mathcal{V}}_n)(\exists\check{W}\in\dot{{\mathcal U}})(V\cap U\neq \emptyset \mbox{ and }U\subseteq \check{W})"$.

Fix $V\in\tau$ nonempty and fix an element $q\in\cohen(\kappa)$. Choose $n$ so that $\vert dom(q)\vert \le n$. Since $\bigcup\mathcal{V}_n$ is dense in $X$, choose $V_{n,x}\in \mathcal{V}_n$ such that $\emptyset\neq V\cap V_{n,x}$. Choose $p\in {\mathbb A}_{n,x}$ such that for an $r\in\cohen(\kappa)$ we have $r\leq q,\, p$. Then 
\[
  r\forces ``\check{V}\cap \check{V}_{n,x}\neq \emptyset \mbox{ and }\check{V}_{n,x}\subseteq \check{V}_{p,x}\in \dot{\mathcal U}"
\]
Thus the set of $r\in\cohen(\kappa)$ forcing the statement in the claim is dense in $\cohen(\kappa)$. The claim follows.

Now let $G$ be a $\cohen(\kappa)$-generic filter and let $V$ be a ground model open set. In the generic extension we find an $n$ and a $U\in\mathcal{V}_n$ and a $W\in \dot{\mathcal U}_G$ with $V\cap U\neq \emptyset$ and $U\subseteq W$. Thus, choose for each $n$ and each $U\in\mathcal{V}_n$ for which this is possible, a $W_U\in\dot{\mathcal U}_G$ with $U\subseteq W_U$. Then $\{W_U: (\exists n)(U\in\mathcal{V}_n)\}$ is a countable subset of $\dot{\mathcal{U}}_G$ and has dense union in $X$. $\Box$

\begin{center}{\bf Preserving almost Lindel\"of}\end{center}

Since ${\sf T}_3$ almost Lindel\"of spaces are Lindel\"of, and since the Lindel\"of property of the product space $^{\omega_1}2$ is ${\sf T}_3$ is not preserved by forcing with countably closed partially ordered sets, the almost Lindel\"of property is not preserved by countably closed forcing. This phenomenon can also occur in generic extensions obtained by forcing with certain countable chain condition partially ordered sets. However, generic extensions by Cohen reals preserve the property of being almost Lindel\"of. 

\begin{lemma}\label{CohenrealClosure}
Let $\kappa>\aleph_0$ be a cardinal number and let $X$ be a topological space and let $U\subseteq  X$ be an open subset of $X$. Then
\[
  {\mathbf 1}_{\cohen(\kappa)}\forces ``\overline{\check{U}} = \check{\overline{U}}".
\]
\end{lemma}
$\Proof$ We show that 
\begin{equation}\label{leftright}
  {\mathbf 1}_{\cohen(\kappa)}\forces ``\overline{\check{U}} \subseteq \check{\overline{U}}".
\end{equation}
and that 
\begin{equation}\label{rightleft}
  {\mathbf 1}_{\cohen(\kappa)}\forces ``\overline{\check{U}} \supseteq \check{\overline{U}}".
\end{equation}

Towards proving (\ref{leftright}), let $\dot{x}$ be a $\cohen(\kappa)$ name for an element of $X$ such that 
${\mathbf 1}_{\cohen(\kappa)}\forces ``\dot{x}\in\overline{\check{U}}"$. Let $A_{\dot{x}}\subseteq\cohen(\kappa)$ be an antichain of $p\in \cohen(\kappa)$ such that there is a corresponding $y_p$, an element of $X$, such that $p\forces``\dot{x}=\check{y}_p"$, and $A_{\dot{x}}$ is maximal with respect to this property.

For any $q\in A_{\dot{x}}$ $y_q$ is a member of $\overline{U}$ since: For each neighborhood $V$ of $y_q$, $q\forces ``\check{V} \mbox{ is a neighborhood of }\dot{x}"$. As $q\forces``\dot{x}\in\overline{\check{U}}"$, we have that $q\forces``\check{V}\cap\check{U}\neq\emptyset"$. Since the parameters in the last sentence forced are all ground model sets, this sentence is true in the ground model.

Suppose that ${\mathbf 1}_{\cohen(\kappa)}$ does not force that $\dot{x}$ is a member of $\check{\overline{U}}$. Choose an $r$ that forces $``\dot{x}\not\in\check{\overline{U}}"$. Choose a $t<r$ and an $s\in X$ such that 
$t\forces``\dot{x}=\check{s}"$. Then in particular we have $s\not\in\overline{U}$. Since $A_{\dot{x}}$ is a maximal antichain of $\cohen(\kappa)$ there is a $q\in A_{\dot{x}}$ which is compatible with $t$: Fix an element $p\in\cohen(\kappa)$ such that $p<q$ and $p<t$. Then $p\forces ``\check{y}_q = \dot{x}$", and it follows that $y_q = s$. But then we find the contradiction that $y_q\in\overline{U}$, while $s\not\in\overline{U}$. This completes the proof of (\ref{leftright}).

Towards proving (\ref{rightleft}), let $\dot{x}$ be a $\cohen(\kappa)$ name for an element of $X$ such that 
${\mathbf 1}_{\cohen(\kappa)}\forces ``\dot{x}\in\check{\overline{U}}"$. Let $A_{\dot{x}}\subseteq\cohen(\kappa)$ be an antichain of $p\in \cohen(\kappa)$ such that there is a corresponding $y_p$, an element of $X$, such that $p\forces``\dot{x}=\check{y}_p"$, and $A_{\dot{x}}$ is maximal with respect to this property. Since for each $p\in A_{\dot{x}}$ we have $p\forces ``\check{y}_p\in \check{\overline{U}}"$ we see that $y_p$ is a member of $\overline{U}$. 

Suppose that ${\mathbf 1}_{\cohen(\kappa)}$ does not force that $\dot{x}\in\overline{\check{U}}$. Let $r$ be such that $r\forces``\dot{x}\not\in\overline{\check{U}}$."  
Let $\dot{V}$ be a $\cohen(\kappa)$ name for a neighborhood of $\dot{x}$ such that $r\forces ``\dot{V}\cap\check{U} = \emptyset"$. 
Since $A_{\dot{x}}$ is a maximal antichain of $\cohen(\kappa)$, choose a $q\in A_{\dot{x}}$ which is compatible with $r$. Then choose $p\in\cohen(\kappa)$ with $p<q$ and $p<r$. Then we have 
\[
  p \forces ``\check{y}_q = \dot{x} \mbox{ and } \dot{V} \mbox{ is a neighborhood of }\dot{x} \mbox{ and }\dot{V}\cap\check{U} = \emptyset."
\]
Choose $t<p$ and an open set $W$ such that $t\forces``\dot{V} =\check{W}"$. Then $t\forces``\check{y}_q\in\check{W}$, and so as $y_q$ is in $\overline{U}$, we have in the ground model that $W\cap U\neq \emptyset$. But then $t\forces ``\dot{V}\cap \check{U}\neq\emptyset"$, a contradiction.
$\Box$

\begin{theorem}\label{cohenrealalmostLpreservation} Let $\kappa>\aleph_0$ be a cardinal. If a topological space $X$ is almost Lindel\"of, then in ${\sf V}^{\cohen(\kappa)}$ $X$ is almost Lindel\"of.  
\end{theorem}
$\Proof$
   For each $n$ let ${\mathcal L}_n$ be an $n$-dowment for $\cohen(\kappa)$
   Let $\tau$ denote the (ground model) topology of $X$. Let $\dot{\mathcal{U}}$ be a $\cohen(\kappa)$ name such that
\[
  {\mathbf {1}}_{\cohen(\kappa)}\forces``\dot{\mathcal{U}}\subseteq \check{\tau} \mbox{ is an open cover of }\check{X}".
\]
  For each element $x$ of $X$ choose a maximal antichain ${\mathbb A}_x\subseteq \cohen(\kappa)$ and ground model open sets $V_{p,x}$, $p\in{\mathbb A}_x$, such that $x\in  V_{p,x}$ and
\[
  p \forces ``\check{V}_{p,x}\in \dot{\mathcal{U}}",
\]
For each $n$ choose ${\mathbb A}_{n,x}\in{\mathcal{L}}_n$ with ${\mathbb A}_{n,x}\subseteq{\mathbb A}_x$ and define $V_{n,x} = \cap\{V_{p,x}: p\in {\mathbb A}_{n,x}\}$. Then $V_{n,x}$ is open, and $\mathcal{U}_n = \{V_{n,x}:x\in X\}$ is in the ground model an open cover of $X$. Since $X$ is almost Lindel\"of, choose a countable set $\mathcal{V}_n\subseteq \mathcal{U}_n$ such that $\{\overline{U}:\, U\in {\mathcal V}_n\}$ is a cover of $X$. 

{\flushleft{\bf Claim: }} ${\mathbf 1}_{\cohen(\kappa)}\forces ``(\forall x\in\check{X})(\exists n)(\exists U\in\mathcal{V}_n)(\exists\check{W}\in\dot{{\mathcal U}})(x\in \overline{U} \mbox{ and }U\subseteq \check{W})"$.

Fix $x\in X$ and fix an element $q\in\cohen(\kappa)$. Choose $n$ so that $\vert dom(q)\vert \le n$. Since $\bigcup\{\overline{U}:\, U\in \mathcal{V}_n\} = X$, choose $V_{n,x}\in \mathcal{V}_n$ such that $x$ is an element of $\overline{V}_{n,x}$. Choose $p\in {\mathbb A}_{n,x}$ such that for an $r\in\cohen(\kappa)$ we have $r\leq q,\, p$. Then 
\[
  r\forces ``x\in \check{\overline{V}}_{n,x} \mbox{ and }\check{V}_{n,x}\subseteq \check{V}_{p,x}\in \dot{\mathcal U}."
\]
Thus by Lemma \ref{CohenrealClosure} the set of $r\in\cohen(\kappa)$ forcing the statement in the claim is dense in $\cohen(\kappa)$. The claim follows.

Now let $G$ be a $\cohen(\kappa)$-generic filter and let $x$ be an element of $X$. In the generic extension we find an $n$ and a $U\in\mathcal{V}_n$ and a $W\in \dot{\mathcal U}_G$ with $x\in \overline{U}$ and $U\subseteq W$. Thus, choose for each $n$ and each $U\in\mathcal{V}_n$ for which this is possible, a $W_U\in\dot{\mathcal U}_G$ with $U\subseteq W_U$. Then $\mathcal{W} =\{W_U: (\exists n)(U\in\mathcal{V}_n)\}$ is a countable subset of $\dot{\mathcal{U}}_G$ and $X = \cup\{\overline{W}:W\in\mathcal{W}\}$. $\Box$

\section{The Rothberger property and weakenings}

For the properties considered in this section the following table contains the property name, its symbolic definition and, where available, a reference for where the property was introduced:
\begin{center}
\begin{tabular}{|l|l|l|}\hline
Property name        & Definition             & Source                 \\ \hline
Rothberger           & $\sone(\open,\open)$   & \cite{R}               \\
Almost Rothberger    & $\sone(\open,\almost)$ & \cite{MSLusin}, p. 251 \\
Weakly Rothberger    & $\sone(\open,\dense)$  & \cite{Daniels}, p. 98  \\ \hline
\end{tabular}
\end{center}
The implications 
\begin{center} Rothberger $\Rightarrow$ Almost Rothberger $\Rightarrow$ weakly Rothberger
\end{center}
are irreversible. It is clear that the Rothberger property implies the Lindel\"of property, the almost Rothberger property implies the almost Lindel\"of property, and the weakly Rothberger property implies the weakly Lindel\"of property.
The weakly Rothberger property and the almost Rothberger property have not been as extensively studied as the Rothberger property. The related properties $\sone(\dense,\dense)$, $\sone(\almost,\dense)$ and $\sone(\almost,\almost)$ have received some attention, but we will not report on these properties here. 

The symbol $\gone^{\omega}(\mathcal{A},\mathcal{B})$ denotes the following game: Players ONE and TWO play $\omega$ innings: In inning $n\in\omega$ ONE first selects an $O_{n}$ from $\mathcal{A}$, and then TWO responds with a $T_{n}\in O_{n}$. A play
\[
  O_0,\, T_0,\, \cdots, O_{n},\, T_{n},\, \cdots
\]
is won by TWO if $\{T_{n}:\, n\in\omega\}\in\mathcal{B}$. Otherwise, ONE wins.

For a specific instance of these games either TWO has a winning strategy in the game $\gone^{\omega}(\mathcal{A},\mathcal{B})$, or ONE has a winning strategy in the game, or else neither player has a winning strategy in the game. In the last case the game is said to be \emph{undetermined}. 

The game $\gone^{\omega}(\open,\open)$ was introduced by Galvin \cite{Galvinpog}. 
In \cite{Galvinpog} Galvin considers winning strategies for TWO in the game $\gone^{\omega}(\open,\open)$ and he proves:
\begin{theorem}[Galvin]\label{galvinth} Let $X$ be a Lindel\"of space such that each one-element subset of $X$ is a ${\sf G}_{\delta}$ subset of $X$. Then the following are equivalent:
   \begin{enumerate}
     \item{TWO has a winning strategy in the game $\gone^{\omega}(\open,\open)$ on $X$.} 
     \item{$X$ is a countable set.} 
   \end{enumerate}
\end{theorem}

In \cite{pawlikowski} Pawlikowski proved the following fundamental theorem, characterizing the Rothberger property in terms of games:
\begin{theorem}[Pawlikowski]\label{pawlikowskith} For a space $X$ the following are equivalent:
\begin{enumerate}
  \item{$X$ has property $\sone(\open,\open)$.} 
  \item{ONE has no winning strategy in the game $\gone^{\omega}(\open,\open)$.} 
\end{enumerate}
\end{theorem}
Thus, in the class of Rothberger spaces whose singleton subsets are ${\sf G}_{\delta}$ sets, the game $\gone^{\omega}(\open,\open)$ is undetermined if, and only if, the space is uncountable. 

Tkachuk considered the game of length $\omega$ for the weakly Rothberger property in \cite{Tkachuk}, denoted there as the game $\Theta_*$. Games seem to be unexamined for the almost Rothberger property. We now explore to what extent analogues of Theorem \ref{galvinth} and Theorem \ref{pawlikowskith} hold for the weakly Rothberger spaces and the almost Rothberger spaces.

\begin{center}{\bf Games and the weakly Rothberger property}\end{center}

The following observation is essentially Proposition 2.6(iii) of \cite{Tkachuk}: 
\begin{lemma}\label{easywinfortwo} If $X$ has a dense subspace $Y$ and TWO has a winning strategy in the game $\gone^{\omega}(\open,\dense)$ on $Y$, then TWO has a winning strategy in $\gone^{\omega}(\open,\dense)$ on $X$. 
\end{lemma}
It follows that on each separable space TWO has a winning strategy in $\gone^{\omega}(\open,\dense)$. We do not currently have a general characterization of spaces for which TWO has a winning strategy in the game $\gone^{\omega}(\open,\dense)$. Theorem \ref{twowinsweaklyrothb} below is the most general result we know that for weakly Rothberger spaces is an analogue of Theorem \ref{galvinth} for Rothberger spaces is . Theorem \ref{twowinsweaklyrothb} has been obtained earlier by Tkachuk and can be deduced from \cite{Tkachuk} Theorem 2.11(i) plus Theorem 3.3(2). Since Tkachuk derives this for a game dual of $\gone^{\omega}(\open,\dense)$ (Theorem 2.11(i) and then introduces and proves the duality (Theorem 3.3(2)), we decided to include a new proof customized to our context here. Also note that although the blanket assumption in \cite{Tkachuk} is that spaces are at least Tychonoff, the arguments in \cite{Tkachuk} also work for the wider class of ${\sf T}_3$ spaces.  
\begin{theorem}[Tkachuk]\label{twowinsweaklyrothb} If $X$ is a first-countable ${\sf T}_3$ space the following are equivalent:
\begin{enumerate}
  \item{$X$ is separable.} 
  \item{TWO has a winning strategy in the game $\gone^{\omega}(\open,\dense)$.} 
\end{enumerate}
\end{theorem}
\Proof The implication (1)$\Rightarrow$(2) follows from the preceding remarks. We must prove that (2)$\Rightarrow(1)$. Thus, let $X$ be a first countable topological space in which TWO has a winning strategy, $\sigma$, in the game $\gone^{\omega}(\open,\dense)$. Let $\prec$ be a well-order of $X$. Let $\theta$ be an infinite regular cardinal which is so large that $X$, its topology, $\open$, $^{<\omega}\open$, $\dense$, $\prec$, $\sigma$, and for each $x\in X$ a neighborhood base $(U_n(x):n\in{\mathbb N})$ of $\{x\}$, are elements of $\htheta$. We may assume that $\mathcal{N} = \{(U_n(x):n\in{\mathbb N}): x\in X\}$ is a member of $\htheta$. Let $(\modelM,\in_{\modelM})$ be a countable elementary submodel of $(\htheta,\in_{\theta})$ such that each of the objects above is an element of $\modelM$. 

{\flushleft{\bf Claim 1: }} For each finite sequence $(\mathcal{U}_1,\,\cdots,\, \mathcal{U}_k)$ of open covers of $X$ there is a point $x\in X$ such that for each neighborhood $U$ of $x$ there is an open cover $\mathcal{U}$ such that $U=\sigma(\mathcal{U}_1,\cdots,\mathcal{U}_k,\mathcal{U})$. 

For suppose the contrary. Then choose for each $x\in X$ a neighborhood $U_x$ such that for each open cover $\mathcal{U}$ of $X$, $U_x\neq \sigma(\mathcal{U}_1,\,\cdots,\, \mathcal{U}_k,\mathcal{U})$. But then as $\mathcal{V} = \{U_x:x\in X\}$ is an open cover of $X$, we have that for some $y\in X$, $\sigma(\mathcal{U}_1,\cdots,\mathcal{U}_k,\mathcal{V})=U_y$. This contradicts the selection of $U_y$, completing the proof of Claim 1.

The sentence in Claim 1 holds in $(\htheta,\in_{\theta})$ and all its parameters are in $\modelM$, so the statement holds in $(\modelM,\in_{\modelM})$. Thus, choose for each finite sequence $\nu$ of open covers of $X$ an $x_{\nu}$ the $\prec$-least element of $X$ satisfying Claim 1. Note that $x_{\nu}$ also satisfies Claim 1 in $(\htheta,\in_{\theta})$, and is in $\htheta$ also the $\prec$-first such element. The set $D = \{x_{\nu}:\nu\in\,^{<\omega}\open\}$ is definable from $\prec$, $\sigma$, $X$ and $^{<\omega}\open$, all parameters in $\modelM$, and thus is a member of $\modelM$.

Enumerate $\open\cap\modelM$ bijectively as $(O_n:n\in{\mathbb N})$. Now $D\cap \modelM$ is the countable set $\{x_{\nu}:\nu\in\, ^{<\omega}\{O_0,O_1,\cdots,O_n,\cdots\}\}$. 
{\flushleft{\bf Claim 2: }} $D\cap \modelM$ is dense in $X$.

For suppose the contrary, and choose a nonempty open set $U$ for which $\overline{U}\cap\modelM\cap D$ is the empty set. The latter is possible since $X$ is ${\sf T}_3$.
 
Since $x_{\emptyset}$ is a member of $\modelM$, and $\mathcal{N}=\{(U_n(x):n\in{\mathbb N}):x\in X\}$ is a member of $\modelM$, also $(U_n(x_{\emptyset}):n\in{\mathbb N})$ (which is definable from $\mathcal{N}$ and $x_{\emptyset}$, both parameters in $\modelM$) is in $\modelM$. But then $\{U_n(x_{\emptyset}):n\in \naturals\}$ is an element, and subset of, $\modelM$. Since $x_{\emptyset}$ is not in $\overline{U}$, fix an $m_1$ with $U_{m_1}(x_{\emptyset})\cap\overline{U}=\emptyset$. Then by Claim 1 $\modelM$ witnesses that there is an open cover $\mathcal{U}$ of $X$ such that $U_{m_1}(x_{\emptyset}) = \sigma(\mathcal{U})$. Thus, choose $n_1$ so that $U_{m_1}(x_{\emptyset}) = \sigma(O_{n_1})$. 

Next, consider $x_{n_1}$. By the same considerations as above there is a neighborhood $U_{m_2}(x_{n_1})$ disjoint from $\overline{U}$, and an $n_2$ such that for the open cover $O_{n_2}$ of $X$, $U_{n_2}(x_{n_1}) = \sigma(O_{n_1},O_{n_2})$. Then apply these considerations to $x_{n_1,n_2}$ to choose a neighborhood $U_{m_3}(x_{n_1,n_2})$ disjoint from $\overline{U}$, and an open cover $O_{n_3}$ of $X$ so that $U_{m_3}(x_{n_1,n_2}) = \sigma(O_{n_1},O_{n_2},O_{n_3})$, and so on. Proceeding like this we obtain a $\sigma$-play
\[
  O_{n_1},\sigma(O_{n_1}),O_{n_2}, \sigma(O_{n_1}, O_{n_2}), \cdots, O_{n_k},\sigma(O_{n_1},\cdots,O_{n_k}),\cdots
\]
for which $\overline{U}\cap \cup_{k\in\naturals}\sigma(O_{n_1},\cdots,O_{n_k})$ is the empty set. Since $\overline{U}$ has nonempty interior this contradicts the fact that $\sigma$ is a winning strategy of TWO in $\gone^{\omega}(\open,\dense)$. 

This completes the proof of $(2)\Rightarrow (1)$.
$\Box$

Next we explore the possibility of an analogue of Theorem \ref{pawlikowskith} for the weakly Rothberger property. In a number of specific examples it has been proven that the property $\sone(\open,\dense)$ is equivalent to ONE not having a winning strategy in the game $\gone^{\omega}(\open,\dense)$. This raises the question of when the property $\sone(\open,\dense)$ is equivalent to player ONE not having a winning strategy in the game $\gone^{\omega}(\open,\dense)$. Here is a partial result in this direction. The Menger property, studied later in this paper, is needed in Theorem \ref{mengerweaklyrothbergerTh} below. The fact we need in the proof of Theorem \ref{mengerweaklyrothbergerTh}, due to Hurewicz, is that a space has the Menger property if, and only if, ONE does not have a winning strategy in the game $\gfin^{\omega}(\open,\open)$. This game is played as follows: ONE and TWO play an inning per $n<\omega$. In the $n$-th inning ONE first chooses an open cover $O_n$ of the space, and then TWO responds by choosing a finite subset $T_n$ of $O_n$. A play $O_0,\, T_0,\, \cdots,\, O_m,\, T_n,\,\cdots$ is won by TWO if $\bigcup_{n<\omega}T_n$ is an open cover of the space. Else, ONE wins.
\begin{theorem}\label{mengerweaklyrothbergerTh}\rm Let $X$ be a Menger space. The following are equivalent:
\begin{enumerate}
\item[(1)] $X$ is weakly Rothberger.
\item[(2)] ONE has no winning strategy in $G^{\omega}_1({\mathcal O},{\mathcal D})$.
\end{enumerate}
\end{theorem}
\Proof We must show that if a space has property $S_1({\mathcal O}, {\mathcal D})$, the ONE has no winning
strategy in the game $G^{\omega}_1({\mathcal O}, {\mathcal D})$. The argument used here is due to Pawlikowski \cite{pawlikowski} for Rothberger spaces. We give some of the details for the convenience of the reader.

Since $X$ is assumed to be a Menger space, it is a Lindel\"{o}f space. Let $F$ be a strategy for ONE in the game $G^{\omega}_1({\mathcal O}, {\mathcal D})$. We may assume that in each inning $F$ calls on ONE to play a countable element of $\mathcal O$.
Define the array $U_{\sigma}$, $\sigma\in \,^{<\omega}{\mathbb N}$, as follows: $(U_n: n\in{\mathbb N})$ enumerates ONE's first move, $F(\emptyset)$. For $n_1$, $(U_{(n_1, n)}: n\in {\mathbb N})$ enumerates $F(U_{n_1})$. For $n_1, n_2$, $(U_{(n_1,n_2,n)}: n\in{\mathbb N})$ enumerates $F(U_{n_1},\,U_{(n_1,n_2)})$, and so on. This array has the property that for each $\sigma$ the set $\{U_{\sigma\frown n}: n\in{\mathbb N}\}$ is in $\mathcal O$.

For fixed $m$ and $j\in{\mathbb N}$ and $\rho$ a function from $\{1,..., j^m\}$ to $\mathbb N$, define the set
$$ U_{\rho}(m, j)=\bigcup_{\sigma\in ^m\{1,...,j\}}(\bigcup\{U_{\sigma \frown\rho\lceil i}:i\leq j^m\})$$
and then for fixed $m$ and $j$ define
\[
  {\mathcal U}(m, j)=\{U_{\rho}(m,j):\rho \mbox{ a function from }\{1,..., j^m\} \mbox{ to }{\mathbb N}\}.
\]
Then each ${\mathcal U}(m, j)\in\mathcal O$.

\bigskip

{\bf Claim:} There exist increasing sequences $\{j_n: n\in{\mathbb  N}\}$ and $\{m_n: n \in {\mathbb N}\}$ such that for each
$x\in X$ and for each $n$ there is a function $\sigma$ from $\{1,...,m_{n+1}-m_n\}$ to
$j_{n+1}$ for which $x\in U_{\sigma}(m_n, j_n)$.

\bigskip

\emph{Proof of the claim.} We observe that $X$ is Menger and this implies that ONE does not have a winning strategy in the corresponding game $G^{\omega}_{fin}({\mathcal O},{\mathcal O})$ using the following strategy, $G$.

For a first move ONE puts $j_1= m_1= 1$, and plays $G(\emptyset)={\mathcal U} (m_1, j_1)$.
For a response $T_1\subseteq {\mathcal U}(m_1, j_1)$ by TWO, ONE first does the following computations: $m_2= m_1+ j_1^{m_1}$, and $j_2 > j_1$ is at least the maximum of all values
of $\sigma$'s for which $U_{\sigma}(m_1, j_1)$ is in $T_1$.
Then ONE plays $G(T_1)={\mathcal U}(m_2,j_2)$.

For a response $T_2\subseteq G(T_1)$ by TWO, ONE again first computes the numbers $m_3$ and $j_3$ according to the
rules that $m_3= m_2 + j^{m_2}_2$, and $j_3 > j_2$ is at least the maximum of all values of $\sigma$'s for
which $U_{\sigma}(m_2, j_2)$ is in $T_2$, and so on.

Look at a $G$-play $G(\emptyset), T_1, G(T_1), T_2, G(T_1,T_2),...$ which is lost by ONE. Then
$\bigcup_{n\in {\mathbb N}} T_n\in \mathcal O$, and we find increasing sequences
$(j_n: n \in {\mathbb  N})$ and $(m_n: n\in {\mathbb  N})$ such that for each $n$:
\begin{enumerate}
\item[(1)] $m_{n+1}=m_n+ j^{m_n}_n$;
\item[(2)] $G(T_1,...T_n)={\mathcal U}(m_{n+1}, j_{n+1})$;
\item[(3)] $j_{n+1}$ is at least as large as the value of an $\sigma$ for which $U_{\sigma}(m_n, j_n)$ is in $T_n$.
\end{enumerate}
It follows that the $m_n$'s and $j_n$'s have the required properties, completing the proof of the claim.

With the sequences $(j_n: n\in {\mathbb N})$ and $(m_n: n \in{\mathbb N})$ fixed, define next for each $n$ the family
${\mathcal W}_n$ as follows: For every sequence $k_1 <...< k_n$ from $\mathbb N$, and for any $\sigma_1,..., \sigma_n$ where each $\sigma_i$ is an $\{1,..., j_{k_i+1}\}$-valued function with domain $m_{k_i+1} - m_{k_i}$, define
\[
  {\mathcal W}(k_1,..., k_n, \sigma_1,..., \sigma_n)=\bigcap_{i\leq n} U_{\sigma_i}(m_{k_i}, j_{k_i}).
\]

${\mathcal W}_n$ consists of all sets of the form $W(k_1,..., k_n, \sigma_1,..., \sigma_n)$.

Since each ${\mathcal W}_n$ is in $\mathcal O$, the selection hypothesis $S_1({\mathcal O}, {\mathcal D})$ applied to $({\mathcal W}_n: n\in{\mathbb N})$ gives for each $n$ a set $S_n=W(k^n_1,...,k^n_n, \sigma^n_1,...,\sigma^n_n)$
such that $\{S_n: n \in {\mathbb N}\}$ is in $\dense$.

Recursively choose for each $n$ an $\ell_n\in\{k^n_1,..., k^n_n\}\backslash \{\ell_i : i < n\}$.
For each $n$ define $\rho_n=\sigma^n_{i_n}$ where $i_n$ is such that $\ell_n= k^n_{i_n}$.

From the definitions we see that for each $n$, $S_n\subseteq U_{\rho_n}(m_{\ell_n},j_{\ell_n})$.

If we now define $f:{\mathbb N}\rightarrow{\mathbb N}$ so that, for each $n$, $f(m_{\ell_n}+i)=\rho_{n}(i)$,
whenever $i\leq m_{\ell_n+1}- m_{\ell_n}$, we find that the play
\[
F(\emptyset), U_{f(1)}, F(U_{f(1)}), U_{f(1),f(2)},F(U_{f(1)}, U_{f(1),f(2)})...
\]
is won by TWO. $\Box$

To what extent is the hypothesis that a space is a Menger space needed in Theorem \ref{mengerweaklyrothbergerTh}? Towards an answer to this question we show that in some models of set theory there are many non-Lindel\"of (and thus non-Menger) spaces where ONE does not have a winning strategy in $\gone^{\omega}(\open,\dense)$. Our proof for one of the main ingredients of the argument is modeled on the arguments in \cite{MStightness}. 

\begin{theorem}\label{cohenrealconversion} Let $\kappa>\aleph_0$ be a cardinal. If topological space $X$ is weakly Lindel\"of, then in ${\sf V}^{\cohen(\kappa)}$ ONE has no winning strategy in $\gone^{\omega}(\open,\dense)$ on $X$. 
\end{theorem}
$\Proof$
Let $\dot{\sigma}$ be a $\cohen(\kappa)$ name such that ${\mathbf 1}_{\cohen(\kappa)}\forces ``\dot{\sigma}$ is a strategy for ONE in $\gone^{\omega}(\open,\dense).$"
By Theorem \ref{cohenrealpreservation},
\[
  {\mathbf 1}_{\cohen(\kappa)}\forces ``\dot{\sigma}(\emptyset) \mbox{ has a countable subset with union dense in }\check{X}."
\]
Choose a $\cohen(\kappa)$ name $\dot{\mathcal{U}}_{\emptyset}$ such that 
\[
  {\mathbf 1}_{\cohen(\kappa)}\forces ``\dot{\mathcal{U}}_{\emptyset}\subseteq \dot{\sigma}(\emptyset) \mbox{ is a countable subset with union dense in }\check{X}."
\]
Thus choose $\cohen(\kappa)$ names $\dot{U}_n$, $n<\omega$ such that 
\[
  {\mathbf 1}_{\cohen(\kappa)}\forces ``\dot{\mathcal{U}}_{\emptyset}=\{\dot{U}_n:n<\omega\}."
\]
Then we have 
\[
  {\mathbf 1}_{\cohen(\kappa)}\forces ``(\forall n)(\dot{\sigma}(\dot{U}_n) \mbox{ has a countable subset with union dense in }\check{X})."
\]
For each $n$ we choose $\cohen(\kappa)$ names $\dot{\mathcal{U}}_n$ and $\dot{U}_{n,k}$, $k<\omega$ such that 
\[
  {\mathbf 1}_{\cohen(\kappa)}\forces ``\dot{\mathcal{U}}_{n}\subseteq \dot{\sigma}(\dot{U}_n) \mbox{ is a countable subset with union dense in }\check{X}"
\]
and
\[
  {\mathbf 1}_{\cohen(\kappa)}\forces ``\dot{\mathcal{U}}_{n}=\{\dot{U}_{n,k}:k<\omega\}."
\]
and so on.
In this way we find for each finite sequence in $\omega$ $\cohen(\kappa)$ names $\dot{\mathcal{U}}_{n_1,\cdots,n_k}$ and $\dot{U}_{n_1,\cdots,n_k}$ such that 
\[
  {\mathbf 1}_{\cohen(\kappa)}\forces ``\{\dot{U}_{n_1,\cdots,n_k,m}:m<\omega\}=\dot{\mathcal{U}}_{n_1,\cdots,n_k}"
\]
and
\[
  {\mathbf 1}_{\cohen(\kappa)}\forces ``\dot{\mathcal{U}}_{n_1,\cdots,n_k}\subseteq \dot{\sigma}(\dot{U}_{n_1},\cdots,\dot{U}_{n_1,\cdots,n_k})"
\]
and
\[
  {\mathbf 1}_{\cohen(\kappa)}\forces ``\dot{\mathcal{U}}_{n_1,\cdots,n_k} \mbox{ is a countable subset with union dense in }\check{X}"
\]

Since $\cohen(\kappa)$ has countable chain condition and each of the names $\dot{\mathcal{U}}_{\tau}$ and $\dot{U}_{\tau}$ is a name for a open single set or a countable set of open sets, there is an $\alpha<\kappa$ such that each of these is a $\cohen(\alpha)$ name. Thus, factoring the forcing as $\cohen(\alpha)*\cohen(\lbrack \alpha,\kappa))$ we may assume that all the named objects are in the ground model. Then, in the generic extension by $\cohen(\lbrack\alpha,\kappa))$ over this ground model there is a function $f\in \,^{\omega}\omega$ such that $f$ is not in any first category set definable from parameters in the ground model.

Now for each nonempty open subset $V$ of $X$ in the ground model define, in the ground model
\[
  F_V = \{f\in\,^{\omega}\omega:(\forall k)(V\cap U_{f\lceil_k}=\emptyset)\}
\]
is first category and is definable from parameters in the ground model only. Thus, in the generic extension by $\cohen(\lbrack \alpha,\kappa))$ $\bigcup_{V\in\check{\tau}}F_V \neq \,^{\omega}\omega$. Choose in this generic extension an $f$ with 
\[
  f\in\,^{\omega}\omega\setminus \bigcup\{F_V: V \mbox{ a ground model open set}\}
\]
Then in the generic extension the $\sigma$-play during which TWO selected the sets $U_{f\lceil_{n}}$, $0<n<\omega$ is won by TWO.
This completes the proof that in the generic extension ONE has no winning strategy in the game $\gone^{\omega}(\open,\dense)$ on $X$.
$\Box$

A weakly Lindel\"of space need not be Lindel\"of: Consider a separable space which is not Lindel\"of. Since every proper forcing preserves being not Lindel\"of, Theorem \ref{cohenrealconversion} shows that it is consistent that there are non-Menger spaces in which ONE has no winning strategy in the game $\gone^{\omega}(\open,\dense)$. It is not clear how much of the Rothberger property must be indirectly present when a space has the weak Rothberger property. Is the converse of the following true?
\begin{theorem}\label{separableRothberger} \rm If $X$ has a dense Rothberger subspace then ONE has no winning strategy in the game $\gone^{\omega}(\open,\dense)$ on $X$.
\end{theorem}

\Proof Let $D$ be a dense Rothberger subspace of $X$. Let $F$ be a strategy for ONE in $\gone^{\omega}(\open,\dense)$ on $X$. From $F$ define a strategy $G$ for ONE on the space $D$ as follows: $G(\emptyset) = \{U\cap D: D\in F(\emptyset)\}$. For $T\in G(\emptyset)$ choose $U_T\in F(\emptyset)$ with $T = D\cap U_T$. Define $G(T) = \{U\cap D:U\in F(U_T)\}$, and so on.

Now by Theorem \ref{pawlikowskith} $G$ is not a winning strategy for ONE in the game $\gone^{\omega}(\open,\open)$. Thus let 
$O_1,\, T_1,\, \cdots,\,O_n,\,T_n,\,\cdots$ 
be a $G$-play of $\gone^{\omega}(\open,\open)$ which is won by TWO on $D$. From the definition of $G$ we find a corresponding sequence $B_1,\, W_1,\, \cdots,\, B_n,\, W_n,\,\cdots$
where $B_1 = F(\emptyset)$ and $O_1 = \{D\cap U:U\in B_1\}$, $W_1\in B_1$ is such that $T_1=D\cap W_1$, and for each $n$ we have $B_{n+1}=F(W_1,\cdots,W_n)$ with $O_{n+1}=\{D\cap U:U\in B_{n+1}\}$ and $T_n = D\cap W_n$. But then we have that 
\[
  D = \bigcup_{n<\infty}T_n\subseteq\bigcup_{n<\infty}W_n.
\]
Since $D$ is dense in $X$, so is $\bigcup_{n<\infty}W_n$. Thus TWO wins this $F$-play of $X$. $\Box$.

\begin{center}{\bf Games and the almost Rothberger property}\end{center}

\begin{theorem}\label{twowinsalmostrothb} Let $(X,\tau)$ be a first countable topological space (but not necessarily ${\sf T}_3$). Then the following are equivalent:
\begin{enumerate}
  \item{$X$ has  countable subset $D$ such that for each neighborhood assignment $f:D\rightarrow\tau$ with $x\in f(x)$ for each $x$ in $D$, the family $\{\overline{f(x)}:x\in D\}$ covers $X$.} 
  \item{TWO has a winning strategy in the game $\gone^{\omega}(\open,\almost)$.} 
\end{enumerate}
\end{theorem}
$\Proof$ (1)$\Rightarrow(2)$: Let a countable set $D$ as in (1) be given, and enumerate it as $(d_n:n<\omega)$. TWO's strategy which chooses in inning $n$ an element $T_n$ of ONE's move $O_n$ so that $d_n\in T_n$ is a winning strategy.
\\ 
Proof of (2)$\Rightarrow(1)$: Let $X$ be a first countable topological space in which TWO has a winning strategy $\sigma$ in the game $\gone(\open,\almost)$. Let $\prec$ be a well-order of $X$. Let $\theta$ be an infinite regular cardinal which is so large that $X$, its topology, $\open$, $^{<\omega}\open$, $\dense$, $\prec$, $\sigma$, and for each $x\in X$ a neighborhood base $(U_n(x):n\in{\mathbb N})$ of $\{x\}$, are elements of $\htheta$. We may assume that $\mathcal{N} = \{(U_n(x):n\in{\mathbb N}): x\in X\}$ is a member of $\htheta$. Let $(\modelM,\in_{\modelM})$ be a countable elementary submodel of $(\htheta,\in_{\theta})$ such that each of the objects above is an element of $\modelM$. 

{\flushleft{\bf Claim 1: }} For each finite sequence $(\mathcal{U}_1,\,\cdots,\, \mathcal{U}_k)$ of open covers of $X$ there is a point $x\in X$ such that for each neighborhood $U$ of $x$ there is an open cover $\mathcal{U}$ such that $U=\sigma(\mathcal{U}_1,\cdots,\mathcal{U}_k,\mathcal{U})$. 

This claim holds because a strategy for TWO in game $\gone^{\omega}(\open,\almost)$ is also a strategy for TWO in the game $\gone^{\omega}(\open,\dense)$. Now use the argument of Claim 1 of Theorem \ref{twowinsweaklyrothb}.

The sentence in Claim 1 holds in $(\htheta,\in_{\theta})$ and all its parameters are in $\modelM$, so the statement holds in $(\modelM,\in_{\modelM})$. Thus, choose for each finite sequence $\nu$ of open covers of $X$ an $x_{\nu}$ the $\prec$-least element of $X$ satisfying Claim 1. Note that $x_{\nu}$ also satisfies Claim 1 in $(\htheta,\in_{\theta})$, and is in $\htheta$ also the $\prec$-first such element. The set $E = \{x_{\nu}:\nu\in\,^{<\omega}\open\}$ is definable from $\prec$, $\sigma$, $X$ and $^{<\omega}\open$, all parameters in $\modelM$, and thus is a member of $\modelM$.

Enumerate $\open\cap\modelM$ bijectively as $(O_n:n\in{\mathbb N})$. Now $D= E\cap \modelM$ is the countable set $\{x_{\nu}:\nu\in\, ^{<\omega}\{O_0,O_1,\cdots,O_n,\cdots\}\}$. 
{\flushleft{\bf Claim 2: }} $D$ has the property defined in (2).

For suppose the contrary. Choose a neighborhood assignment $f:D\rightarrow\tau$ for which $\{\overline{f(x)}:x\in D\}$ does not cover $X$. Pick a point $y\in X\setminus\bigcup\{\overline{f(x)}:x\in D\}$.

Since $x_{\emptyset}$ is a member of $\modelM$, and $\mathcal{N}=\{(U_n(x):n\in{\mathbb N}):x\in X\}$ is a member of $\modelM$, also $(U_n(x_{\emptyset}):n\in{\mathbb N})$ (which is definable from $\mathcal{N}$ and $x_{\emptyset}$, both parameters in $\modelM$) is in $\modelM$. But then $\{U_n(x_{\emptyset}):n\in \naturals\}$ is an element, and subset of, $\modelM$. Fix an $m_1$ with $U_{m_1}(x_{\emptyset})\subset f(x_{\emptyset})$. Since $y$ is not in $\overline{f(x_{\emptyset})}$, $y$ is also not a member of $\overline{U}_{m_1}(x_{\emptyset})$. 

Then by Claim 1 $\modelM$ witnesses that there is an open cover $\mathcal{U}$ of $X$ such that $U_{m_1}(x_{\emptyset}) = \sigma(\mathcal{U})$. Thus, choose $n_1$ so that $U_{m_1}(x_{\emptyset}) = \sigma(O_{n_1})$. 

Next, consider $x_{n_1}$. By the same considerations as above there is a neighborhood $U_{m_2}(x_{n_1})$ with $y\not\in\overline{U_{m_2}(x_{n_1})}$, and an $n_2$ such that for the open cover $O_{n_2}$ of $X$, $U_{m_2}(x_{n_1}) = \sigma(O_{n_1},O_{n_2})$.

Apply these considerations to $x_{n_1,n_2}$ to choose a neighborhood $U_{m_3}(x_{n_1,n_2})$ with $y$ not a member of $\overline{U}_{m_3}(x_{n_1,n_2})$ and an open cover $O_{n_3}$ of $X$ so that $U_{m_3}(x_{n_1,n_2}) = \sigma(O_{n_1},O_{n_2},O_{n_3})$, and so on. Proceeding like this we obtain a $\sigma$-play 
\[
  O_{n_1},\sigma(O_{n_1}),O_{n_2}, \sigma(O_{n_1}, O_{n_2}), \cdots, O_{n_k},\sigma(O_{n_1},\cdots,O_{n_k}),\cdots
\] 
for which $y$ is not a member of $\cup_{k\in\naturals}\overline{\sigma(O_{n_1},\cdots,O_{n_k})}$. This contradicts the fact that $\sigma$ is a winning strategy of TWO in $\gone(\open,\almost)$. 

This completes the proof of $(2)\Rightarrow (1)$. $\Box$

\begin{corollary}\label{T2almostrothb} If $X$ is a first countable ${\sf T}_2$-space such that TWO has a winning strategy in the game $\gone^{\omega}(\open,\almost)$, then $X$ is countable.
\end{corollary}
\Proof: Let $D\subseteq X$ be the countable subset as in (2) of Theorem \ref{twowinsalmostrothb}. We claim that $D=X$. For suppose the contrary and choose $y\in X\setminus D$. Then as $X$ is ${\sf T}_2$, choose for each $d\in D$ a neighborhood $f(d)$ with $y$ not in $\overline{f(d)}$. It follows that $f$ is a neighborhood assignment violating the property (2) of $D$, a contradiction. $\Box$

In Theorem \ref{twowinsalmostrothb} some hypothesis like first countability is needed to derive that (2) implies (1): In Example F in the examples section of the paper we present an uncountable ${\sf T}_2$ space which is not Lindel\"of, not first countable, and TWO has a winning strategy in $\gone^{\omega}(\open,\almost)$, and the space does not satisfy property (1) of Theorem \ref{twowinsalmostrothb}. 

In Theorem 8 of \cite{MSRothbgroup} it is shown that there is for each infinite cardinal number $\kappa$ a ${\sf T}_4$ Lindel\"of P-space $X$ of cardinality $\kappa$ for which TWO has a winning strategy in the game $\gone^{\omega}(\open,\open)$, and thus also in the games $\gone^{\omega}(\open,\dense)$ and $\gone^{\omega}(\open,\almost)$. By Galvin's Theorem, Theorem \ref{galvinth}, $X$ has some one-element subset which is not a ${\sf G}_{\delta}$ set. Since $X$ is in fact a topological group it follows that no one-element subset of $X$ is a ${\sf G}_{\delta}$ set.

\begin{question}\label{almostptsgdelta} Can there be an uncountable ${\sf T}_2$-space $X$ such that each one-element subset of $X$ is a ${\sf G}_{\delta}$ set, and TWO has a winning strategy in the game $\gone^{\omega}(\open,\almost)$?
\end{question}

\begin{question}\label{almostT1} Can there be an uncountable ${\sf T}_1$-space $X$ which is first countable such that TWO has a winning strategy in the game $\gone^{\omega}(\open,\almost)$?
\end{question}

\begin{theorem}\label{mengeralmostrothbergerTh}\rm Let $X$ be a Menger space. The following are equivalent:
\begin{enumerate}
\item[(1)] $X$ is almost Rothberger.
\item[(2)] ONE has no winning strategy in $\gone^{\omega}({\mathcal O},\almost)$.
\end{enumerate}
\end{theorem}
\Proof The proof of (1) implies (2) proceeds like the proof of the corresponding implication in the proof of Theorem \ref{mengerweaklyrothbergerTh}. But at the stage of that proof where $\sone(\open,\dense)$ is applied to the sequence of ${\mathcal W}_n$'s, apply instead $\sone(\open,\almost)$. $\Box$

Using the same technique as in the proof of Theorem \ref{cohenrealconversion} we find:
\begin{theorem}\label{cohenrealalmostconversion} For $X$ an almost Lindel\"of space and $\kappa$ an uncountable cardinal, 
\[
  {\mathbf 1}_{\cohen(\kappa)}\forces``\mbox{ONE has no winning strategy in the game }\gone^{\omega}(\open,\overline{\open}) \mbox{ on }\check{X}."
\]
\end{theorem}

There are non-Lindel\"of, almost Lindel\"of spaces. Since Cohen real forcing preserves non-Lindel\"of, Theorem \ref{cohenrealalmostconversion} implies it is consistent that there are non-Menger spaces for which ONE has no winning strategy in the game $\gone^{\omega}(\open,\almost)$.

\section{The Menger property and weakenings}

The following table contains names for properties considered in this section, symbolic definitions, and where available a reference for where the property was introduced:
\begin{center}
\begin{tabular}{|l|l|l|}\hline
Property name        & Definition             & Source                 \\ \hline
Menger               & $\sfin(\open,\open)$   & \cite{Hurewicz}        \\
Almost Menger        & $\sfin(\open,\almost)$ & \cite{Ko2}             \\
Weakly Menger        & $\sfin(\open,\dense)$  & \cite{Daniels}, p. 94  \\ \hline
\end{tabular}
\end{center}
The weakly Menger property was investigated in \cite{P}. Figure \ref{weakrelfig} below depicts the implications among the properties have been introduced thus far:

\setlength{\unitlength}{.5in}
\begin{center}
\begin{figure}[h]
\begin{picture}(6,4.05)(-1,1.75)
\put(0,2.25){$\sone(\open,\open)$}

\put(0.45,2.5){\vector(0,1){1}}

\put(1,2.35){\vector(1,0){1.45}}
\put(2.75,2.25){$\sone(\open,\almost)$}

\put(3.2,2.5){\vector(0,1){1}}

\put(3.85,2.35){\vector(1,0){1.45}}
\put(5.5,2.25){    $\sone(\open,\dense)$}

\put(5.95,2.5){\vector(0,1){1}}
\put(0,3.75){$\sfin(\open,\open)$}

\put(0.45,4.05){\vector(0,1){1}}

\put(1.4,3.85){\vector(1,0){1.05}} 
\put(2.75,3.75){$\sfin(\open,\almost)$}

\put(3.2,4.05){\vector(0,1){1}}

\put(4.15,3.85){\vector(1,0){1.15}}  
\put(5.5, 3.75){$\sfin(\open,\dense)$}

\put(5.95, 4.05){\vector(0,1){1.05}}
\put(0,5.25){Lindel\"of}

\put(1.2,5.3){\vector(1,0){1.25}}

\put(2.75,5.25){Almost Lindel\"of}

\put(4.75,5.3){\vector(1,0){0.45}}


\put(5.5,5.25){Weakly Lindel\"of}

\end{picture}
\caption{Basic Relationships \label{weakrelfig}}
\end{figure}
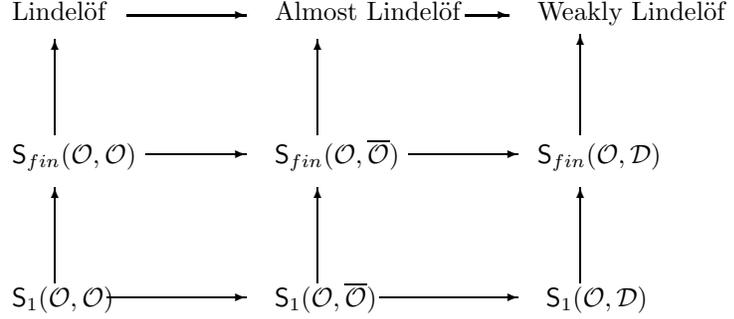
\end{center}

If a space is ${\sf T}_3$, then it is almost Menger if, and only if, it is Menger. However, $\sfin(\open,\almost)$ is not equivalent to $\sfin(\almost,\almost)$, even in the context of metrizable spaces:
\begin{theorem}[\cite{MSLusin}, Theorem 2] An uncountable set of real numbers has the property $\sfin(\almost,\almost)$ if, and only if, it is a Lusin set.
\end{theorem}

\begin{center}{\bf Games and the Menger property}\end{center}

In Corollary 4 of \cite{Telgarsky} Telg\'arsky characterized the metrizable spaces for which TWO has a winning strategy as the ones that are $\sigma$-compact. A more direct proof of Telg\'arsky's result was given in \cite{MST}. There appears to be no satisfactory more general characterizations of topological spaces for which player TWO has a winning strategy in the games $\gfin^{\omega}(\open,\open)$, $\gfin^{\omega}(\open,\almost)$ or $\gfin^{\omega}(\open,\dense)$.  

It is not clear how this characterization would generalize to say ${\sf T}_4$ spaces. It is false that for ${\sf T}_4$-spaces TWO has a winning strategy in the game $\gfin^{\omega}(\open,\open)$ if, and only if, the space is $\sigma$-compact. This can be seen by considering the examples in Section 4 of \cite{MSRothbgroup}: By Theorem 8 of \cite{MSRothbgroup} there is for each infinite cardinal number $\kappa$ a ${\sf T}_4$ Lindel\"of {\sf P}-group of cardinality $\kappa$ such that TWO has a winning strategy in the game $\gone^{\omega}(\open,\open)$ (and thus in $\gfin^{\omega}(\open,\open)$). An uncountable ${\sf T}_4$ Lindel\"of {\sf P}-space cannot be a closed subset of a $\sigma$-compact space. 

It is clear that if a space is $\sigma$-compact, then TWO has a winning strategy in the game $\gfin^{\omega}(\open,\open)$. A weaker condition than this permits TWO a winning strategy in the game $\gfin^{\omega}(\open,\dense)$: 
\begin{definition}\rm \cite{PW}
$X$ is {\sf H}-\emph{closed }if every open cover $\mathcal U$ of $X$ has a finite subfamily ${\mathcal V}$
whose union is dense in $X$ (i.e. $X\subseteq cl_X(\bigcup_{V\in {\mathcal V}}V)$).
\end{definition}
It is well-known that a ${\sf T}_2$-space is {\sf H}-closed if, and only if it is a closed subspace of each ${\sf T}_2$-space it embeds in. Equally well-known, a ${\sf T}_3$-space is {\sf H}-closed if, and only if, it is compact.

\begin{definition}\rm 
$X$ is $\sigma$-{\sf H}-\emph{closed }if it is a union of countably many subspaces, each of which is {\sf H}-closed.
\end{definition}

\begin{proposition}
\rm
If $X$ contains a dense $\sigma$-{\sf H}-closed subspace then TWO has a winning strategy in the game $\gfin^{\omega}(\open,\dense)$.
\end{proposition}

\Proof Let $Y\subseteq X$ be a dense $\sigma$-{\sf H}-closed subspace of $X$. 
Since $Y$ is a $\sigma$-H-closed space, write $Y = \bigcup_{n\in{\mathbb N}}Y_n$ where each $Y_n$ is an {\sf H}-closed subspace of $Y$. Define a strategy $\sigma$ for player TWO as follows: When in inning $n$ ONE plays an open cover $O_n$ of $X$, let $\sigma(O_1,\cdots,O_n)\subseteq O_n$ be a finite set with $\overline{\bigcup \sigma(O_1,\cdots,O_n)}\supseteq Y_n$. Then $\sigma$ is a winning strategy for TWO. $\Box$

\begin{question}($T_2$) Is it true that $X$ contains a dense $\sigma$-{\sf H}-closed subspace if, and only if, TWO has a winning strategy in $G_{fin}^{\omega}({\mathcal O}, {\mathcal D})$?
\end{question}

\begin{question}($T_3$) Is it true that $X$ has a dense $\sigma$-compact subset if, and only if, TWO has a winning strategy in $G_{fin}^{\omega}({\mathcal O}, {\mathcal D})$?
\end{question}

\begin{question}
Characterize the topological spaces for which TWO has  winning strategy in:
\begin{enumerate}
  \item{the game $\gfin^{\omega}(\open,\open)$.} 
  \item{the game $\gfin^{\omega}(\open,\almost)$.}
  \item{the game $\gfin^{\omega}(\open,\dense)$.} 
\end{enumerate}
\end{question}

Conditions under which ONE has no winning strategy in these games are somewhat better understood. 
Hurewicz, who introduced and studied the Menger property in \cite{Hurewicz}, proved there
\begin{theorem}[Hurewicz]\label{Hurewiczplayerone} For a space $X$ the following are equivalent:
\begin{enumerate}
  \item{$X$ has the Menger property} 
  \item{ONE has no winning strategy in the game $\gfin^{\omega}(\open,\open)$.} 
\end{enumerate}
\end{theorem}

For Lindel\"of spaces Hurewicz's proof of Theorem \ref{Hurewiczplayerone} generalizes to also give the corresponding characterizations for weakly Menger and almost Menger spaces. For the convenience of the reader we now give the proof of the characterization of weakly Menger Lindel\"of spaces, and then indicate what modification is needed to also obtain the characterization for Lindel\"of almost Menger spaces.

\begin{theorem}\label{weakmengergame} Let $X$ be a Lindel\"of space. Then the following are equivalent:
\begin{enumerate}
  \item{$X$ is weakly Menger.}
  \item{ONE does not have a winning strategy in the game $\gfin^{\omega}(\open,\dense)$.}   
\end{enumerate}  
\end{theorem}
\Proof The implication that if a Lindel\"of space satisfies $S_{fin}({\mathcal O}, {\mathcal D})$ then ONE has no winning
strategy in the game $\gfin^{\omega}({\mathcal O},{\mathcal D})$ requires proof. The argument used here is due to Hurewicz \cite{Hurewicz} for Menger spaces, and has been used in \cite{COC5} for a different context. We give some of the details of Hurewicz's argument for the convenience of the reader. 

Let $X$ be Lindel\"{o}f and let $F$ be a strategy for ONE. Without loss of generality, we may assume that each move of ONE according to the strategy $F$, is an ascending $\omega$-chain of open sets covering $X$. 

Write $F(\emptyset) = (U(n): n\in {\mathbb N})$, listed in $\subset$-increasing order. Then, for each $n$, list $F(U_{(n)})$ in $\subset$-increasing order as $(U_{(n,m)}: m \in {\mathbb N})$, and so on.

Supposing that $U_{\tau}$, has been defined for each finite sequences $\tau$ of length at most $k$ of positive integers, we now define for each $(n_1 , . . . , n_k)$:

\[
  F(U_{(n_1)}, . . . , U_{(n_1, ,,..., n_k)}) = (U_{(n_1,..., n_k,m)}: m \in \naturals).
\]

Then the family
$(U_{\tau}: \tau$ a finite sequence of positive integers$)$ has the following properties for each $\tau$:
\begin{enumerate}
\item[(1)] If $m$ is less than $n$, then $U_{\tau\frown (m)}$ is a subset of $U_{\tau\frown (n)}$.
\item[(2)] For each $n$, $U_{\tau}\subseteq U_{\tau\frown (n)}$.
\item[(3)] $\{U_{\tau\frown (n)}: n$ a positive integer$\}$ is an open cover of $X$.
\end{enumerate}

Define for each $n$ and $\kappa$:

\[
  U^n_k=\left\{
          \begin{array}{ll}
            U_{(k)}, & \hbox{if n=1;} \\
            (\bigcap \{U_{\tau\frown (k)}:\tau\in ^{n-1}{\mathbb N}\})\cap U^{n-1}_{\kappa}, & \hbox{otherwise.}
          \end{array}
        \right.
\]

Note that for each $n$ the set $\{U^n_{k}:k \in {\mathbb N}\}$, denoted ${\mathcal U}_n$, is an open
cover of $X$. To see this, first show (by induction) that for each $(i_1, . , i_n)$ such that max$\{i_l, . . . , i_n\} \geq k$ one has $U^n_{k}\subseteq U_{(i_1,..., i_n)}$. It then follows that each $U^n_k$ is an intersection of finitely many open sets, and thus is itself open. Now observe that by its definition each ${\mathcal U}_n$ is an increasing chain of open sets such that each ${\mathcal U}_n$ is an open cover of $X$.

Apply the fact that $X$ is weakly Menger to the sequence $({\mathcal U}_n:n\in \naturals)$. Since each ${\mathcal U}_n$ is an ascending chain this gives for each $n$ a $U^n_{k_n}\in {\mathcal U}_n$ such that $\{U^n_{k_n}:n\in \naturals\}$ is in $\mathcal D$. Since, for each $n$, $U^n_{k_n}\subseteq U_{(k_1,...,k_n)}$, the sequence of moves
$U_{k_1},U_{k_1,k_2},\cdots$ by TWO defeats ONE's strategy $F$. $\Box$

The corresponding theorem for almost Menger is:
\begin{theorem}\label{almostmengergame} Let $X$ be a Lindel\"of space. Then the following are equivalent:
\begin{enumerate}
  \item{$X$ is almost Menger.}
  \item{ONE does not have a winning strategy in the game $\gfin^{\omega}(\open,\almost)$.}   
\end{enumerate}  
\end{theorem}
\Proof The proof proceeds like for Theorem \ref{weakmengergame}. In that proof at one point we apply the selection property $\sfin(\open,\dense)$ to a sequence $(\mathcal{U}_n:n<\omega)$ of special open covers of $X$. Applying the selection property $\sfin(\open,\almost)$ instead at the same point of the proof produces a proof of Theorem \ref{almostmengergame}.
$\Box$

The hypothesis that $X$ is Lindel\"of in Theorems \ref{weakmengergame} or \ref{almostmengergame} is not necessary: By Theorem \ref{cohenrealconversion} it is consistent that there is a non-Lindel\"of space $X$ for which TWO has a winning strategy in $\gone^{\omega}(\open,\dense)$, and thus in $\gfin^{\omega}(\open,\dense)$. Similar remarks apply to the almost Menger case. We do not know the answers to the following questions:
\begin{question}
Is there an almost Menger space for which ONE has a winning strategy in the game $\gfin^{\omega}(\open,\almost)$?
\end{question}

\begin{question}
Is there a weakly Menger space for which ONE has a winning strategy in the game $\gfin^{\omega}(\open,\dense)$?
\end{question}

\section{Examples}

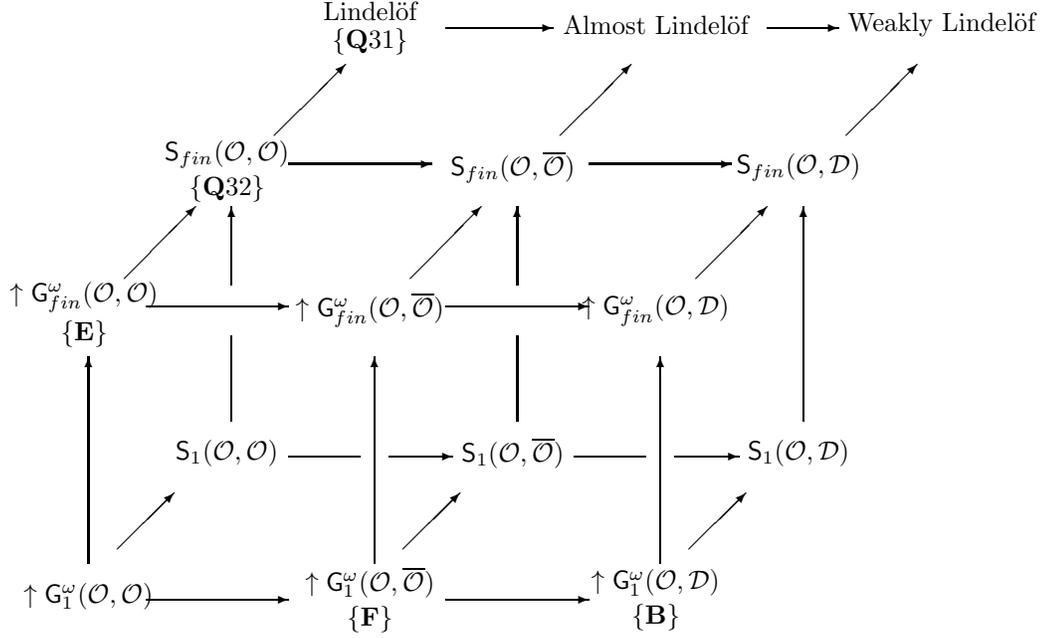
\begin{figure}[h]
\unitlength=.95mm
\begin{picture}(140.00,95.00)(10,10)
\put(20.00,20.00){\makebox(0,0)[cc]
                  {\shortstack {$\uparrow\gone^{\omega}(\open,\open)$\\   {}        }}} 

\put(60.00,20.00){\makebox(0,0)[cc]
{\shortstack {$\uparrow\gone^{\omega}(\open,\almost)$\\ $\{ {\bf F}\} $   } }}

\put(100.00,20.00){\makebox(0,0)[cc]
{\shortstack {$\uparrow\gone^{\omega}(\open,\dense)$\\ $\{ {\bf B} \}$   } }}
\put(20.00,60.00){\makebox(0,0)[cc]
{\shortstack {$\uparrow\gfin^{\omega}(\open,\open)$\\ $\{  {\bf E} \} $   } }}

\put(60.00,60.00){\makebox(0,0)[cc]
{\shortstack {$\uparrow\gfin^{\omega}(\open,\almost)$\\ {}   } }} 

\put(100.00,60.00){\makebox(0,0)[cc]
{\shortstack {$\uparrow\gfin^{\omega}(\open,\dense)$\\ {} } }} 
\put(40.00,40.00){\makebox(0,0)[cc]
{\shortstack {$\sone(\open,\open)$\\  {}  } }} 

\put(80.00,40.00){\makebox(0,0)[cc]
{\shortstack {$\sone(\open,\almost)$\\ {}   } }} 

\put(120.00,40.00){\makebox(0,0)[cc]
{\shortstack {$\sone(\open,\dense)$\\ {}  } }}  
\put(40.00,80.00){\makebox(0,0)[cc]
{\shortstack {$\sfin(\open,\open)$\\ $\{ {\bf Q \ref{Q33}} \}$  } }}

\put(80.00,80.00){\makebox(0,0)[cc]
{\shortstack {$\sfin(\open,\almost)$\\ {}  } }} 

\put(120.00,80.00){\makebox(0,0)[cc]
{\shortstack {$\sfin(\open,\dense)$\\ {}  } }} 
\put(60.00,100.00){\makebox(0,0)[cc]
{\shortstack {Lindel\"of\\ $\{ {\bf Q \ref{Q32}} \}$   } }}

\put(100.00,100.00){\makebox(0,0)[cc]
{\shortstack {Almost Lindel\"of\\ {}  } }} 

\put(140.00,100.00){\makebox(0,0)[cc]
{\shortstack {Weakly Lindel\"of\\ {} } }} 
\put(70.00,100.00){\vector(1,0){15.00}}
\put(115.00,100.00){\vector(1,0){10.00}}
\put(46.00,85.00){\vector(1,1){10.00}}
\put(86.00,85.00){\vector(1,1){10.00}}
\put(126.00,85.00){\vector(1,1){10.00}}
\put(25.00,65.00){\vector(1,1){10.00}}   
\put(105.00,65.00){\vector(1,1){10.00}}  
\put(65.00,65.00){\vector(1,1){10.00}}
\put(28.00,61.00){\vector(1,0){20.00}}
\put(70.00,61.00){\vector(1,0){20.00}}
\put(48.00,81.00){\vector(1,0){20.00}} 
\put(90.00,81.00){\vector(1,0){20.00}} 
\put(40.00,45.00){\line(0,1){12.00}}
\put(80.00,45.00){\line(0,1){12.00}}   
\put(120.00,45.00){\vector(0,1){30.00}}
\put(40.00,64.00){\vector(0,1){11.00}}
\put(80.00,64.00){\vector(0,1){11.00}}
\put(20.00,25.00){\vector(0,1){29.00}}
\put(60.00,25.00){\vector(0,1){29.00}}
\put(24.00,27.00){\vector(1,1){8.00}}
\put(64.00,27.00){\vector(1,1){8.00}}
\put(104.00,27.00){\vector(1,1){8.00}}
\put(100.00,25.00){\vector(0,1){29.00}}
\put(28.00,20.00){\vector(1,0){20.00}}
\put(70.00,20.00){\vector(1,0){20.00}}
\put(48.00,40.00){\line(1,0){10.00}} 
\put(88.00,40.00){\line(1,0){10.00}} 
\put(62.00,40.00){\vector(1,0){9.00}}
\put(102.00,40.00){\vector(1,0){9.00}} 
\end{picture}
\caption{Classes considered in this paper \label{cshl}}
\end{figure}

The relationships among the covering properties we consider in this paper are indicated in the Figure \ref{weakrelfig}. Also considering the game theoretic versions, the strongest property considered in this paper is the property that TWO has a winning strategy in the game $\gone^{\omega}(\open,\open)$. The following diagram, Figure \ref{cshl}, is Figure \ref{weakrelfig} updated to include the classes where TWO has a winning strategy in the corresponding game. We use the symbol $\uparrow\gfin^{\omega}(\mathcal{A},\mathcal{B})$ to denote that TWO has a winning strategy in the corresponding game.


We now consider examples that distinguish these classes from each other. We are missing two examples, as indicated in the following two questions:

\begin{question}\label{Q32} Is there a Lindel\"of space which is not weakly Menger?
\end{question}

\begin{question}\label{Q33} Is there a Menger space for which TWO does not have a winning strategy in $\gfin^{\omega}(\open,\dense)$?
\end{question}
These questions are associated with the corresponding vertices in Figure \ref{cshl}. The rest of our examples are as follows:

{\flushleft{\bf A:}} It is easy to find compact ${\sf T}_2$ spaces of arbitrary large infinite cardinality in this class: For an infinite cardinal number $\kappa$ let ${\sf D}_{\kappa}$ be the one-point compactification of a discrete space of cardinality $\kappa$. Then TWO has a winning strategy in $\gone^{\omega}(\open,\open)$ on ${\sf D}_{\kappa}$.

{\flushleft{\bf B:} Non-almost Lindel\"of spaces where TWO has a winning strategy in $\gone^{\omega}(\open,\dense)$.}\\
Note that such an example indicates that none of the implications from the middle panel to the right panel of Figure \ref{cshl} is reversible. \\
(a) The space ${\mathbb P}$ of irrational numbers with the topology inherited from the real line is not Menger. As ${\mathbb P}$ is ${\sf T}_3$, it follows that this space is also not almost Menger. For any refinement of this topology on ${\mathbb P}$, these statements remain true. As ${\mathbb P}$ separable, TWO has a winning strategy in the game $\gone^{\omega}(\open,\dense)$.
Define the topology
\[
  \tau_0 :=\{U\setminus C: U \mbox{ open in the usual topology on }{\mathbb P} \mbox{ and } C\subseteq {\mathbb P} \mbox{ is countable}\}
\]
Then $({\mathbb P},\tau_0)$ is a ${\sf T}_2$-space but no longer a ${\sf T}_3$-space. It is still Lindel\"of and not almost Menger, but it is no longer separable. Yet, TWO has a winning strategy in the game $\gone^{\omega}(\open,\dense)$ on $({\mathbb P},\tau_0)$. This indicates that in Theorem \ref{twowinsweaklyrothb} the hypothesis that the space is ${\sf T}_3$ cannot be weakened to ${\sf T}_2$.\\
(b) Examination of Examples 6 in \cite{P} shows that also that example is non-separable, ${\sf T}_2$ but not ${\sf T}_3$, and the space is not almost Lindel\"of, but TWO has a winning strategy in the game $\gone^{\omega}(\open,\dense)$. \\ 
(c) Example 11 of \cite{P} is, on the other hand, a separable ${\sf T}_{3\frac{1}{2}}$-space which is not almost Lindel\"of.\\
(d) Let ${\mathbb S}$ denote the \emph{Sorgenfrey line}, the topological space obtained from refining the standard topology on the real line by also declaring each interval of the form $\lbrack a,\, b)$ open. Since the set of rational numbers still is a dense subset of ${\mathbb S}$,  TWO has a winning strategy in $\gone^{\omega}(\open,\dense)$ on ${\mathbb S}$ and all its finite powers. In Lemma 17 of \cite{LBMS} it was shown that ${\mathbb S}$ does not have the property $\sfin(\open,\open)$, and since ${\mathbb S}$ is ${\sf T}_3$, this means that ${\mathbb S}$ is not almost Menger. Note that finite powers of ${\mathbb S}$ are still separable, but are not Lindel\"of and as these powers remain ${\sf T}_3$, they also are not almost Lindel\"of. Thus, for finite powers of ${\mathbb S}$ TWO has a winning strategy in $\gone^{\omega}(\open,\dense)$ while the space is not almost Lindel\"of.\\
(e) The space ${\mathbb Z}^{\omega_1}$ is ${\sf T}_3$ and not ${\sf T}_4$, it is not Lindel\"of and thus also not almost Lindel\"of. But it is separable, and so TWO has a winning strategy in the game $\gone^{\omega}(\open,\dense)$ for this space.

{\flushleft{\bf C:} Non-Lindel\"of (but almost Menger) spaces where TWO has a winning strategy in $\gone^{\omega}(\open,\dense)$.}\\
Such an example shows that none of the implications from the top left edge of Figure \ref{cshl} to the top middle is reversible.\\
Example 77 of \cite{SS}, the deleted radius topology in the plane, is a non-Lindel\"of space. In \cite{K}, Example 1, Ko\v{c}ev points out that the space in this example is almost Menger (and thus not ${\sf T}_3$). Indeed, TWO has a winning strategy in the game $\gfin^{\omega}(\open,\almost)$: When ONE presents TWO with an open cover, TWO may first replace it with basic open sets consisting of the appropriate deleted radius open disks. Ignoring the deleted radii, this would be a move of ONE in the usual topology of $\reals^2$, and TWO may, in the $n$-th inning, choose finitely many of these open disks (including radii) that cover the Euclidean compact set $\lbrack-n,\,n\rbrack\times \lbrack-n,\,n\rbrack$. Then remove the radii so as to recover sets from the replacement family of basic open sets of the deleted radius topology, and then select finitely many elements of ONE's presented cover that contain the corresponding finitely many basis elements. This example illustrates that in the proof of Theorem \ref{almostmengergame} the hypothesis that a space be Lindel\"of is not required for the conclusion that TWO has a winning strategy in the almost Menger game.

This space is not almost Rothberger, as can be shown by for each positive integer $n$  letting $\mathcal{U}_n$ be the open cover consisting of all open discs of area less than $\frac{1}{2^n}$ with horizontal radius (excluding the center) removed. If for each $n$ we choose a $U_n\in\mathcal{U}_n$, then the total area covered by the sets $\overline{U}_n$, $n<\infty$, is finite, and so these sets do not cover the plane. But the set $D$ of points in the plane with rational coordinates only is still dense in the deleted radius topology, so that this space is separable. It follows that TWO has a winning strategy in the game $\gone^{\omega}(\open,\dense)$. Should we refine the deleted radius topology further by declaring all countable sets closed, the resulting space would no longer be separable, but TWO would still have a winning strategy in the game $\gone^{\omega}(\open,\dense)$, and in the game $\gfin^{\omega}(\open,\almost)$.

{\flushleft{\bf D:} Non-Lindel\"of almost Rothberger spaces.}\\
Such an example eliminates one more implication from the middle panel to the left panel of Figure \ref{cshl}. \\
Consider subspaces of the space in Example 77 of \cite{SS}, the Euclidean plane with the deleted radius topology. It is a non-Lindel\"of space. Now assume $X$ is an uncountable subset of $\reals$ such that $X\times X$ has the Rothberger property. Let ${\sf R}(X)$ denote the set $X\times X$ with the deleted radius topology inherited. Then one can show that ${\sf R}(X)$ still is not Lindel\"of, but it is almost Rothberger.

{\flushleft{\bf E:} Compact ${\sf T}_2$ spaces which are not weakly Rothberger.}\\
This example shows that there are no implications form the top level to the bottom level of Figure \ref{cshl}.\\
Let $\bf I$ be the closed unit interval. Let $X$ be a dense subset of $\bf I$. Consider the following subspace $T(X)$ of the \emph{Alexandroff double of} {\bf I}. $T(X) = {\bf I}\times\{0\}\cup  X \times\{l\}$.
For $A\subset {\bf I}$ and for $x\in {\bf I}$ we write $A_i$ for $A \times\{i\}$ and $x_i$ for $(x,i)$, $i\in\{0, 1\}$.
The family ${\mathcal B}=\{U_0 \cup ((U \cap X)_1\backslash\{x_1\}): U$ open in ${\bf I}$ and $x\in U\cap X\} \cup \{\{x_1\}:x\in X\}$
is a basis for a topology on $T(X)$. In this topology $T(X)$ is compact and $T_2$. In particular, TWO has a winning strategy in the game $\gfin^{\omega}(\open,\open)$.

The weakly Rothberger property for $T(X)$ is connected with the classical strong measure zero sets of real numbers.
A subset $X$ of the real line $\mathbb R$ has \emph{strong measure zero} if there is for every sequence $(\varepsilon_n: n\in{\mathbb N})$ of positive real numbers a sequence
$(J_n: n\in{\mathbb N})$ of nonempty open intervals such that each $J_n$ has length at most $\varepsilon_n$, and
$X \subseteq \bigcup_{n\in {\mathbb N}} J_n$. 
This concept was introduced in \cite{Bo} where Bore1 observed that every countable set of real numbers has this property.
Borel conjectured that every strong measure zero set is countable. The truth of tis conjecture is independent of {\sf ZFC}.
In \cite{COC4} it is shown that for X a dense subset of $\bf I$, the following are equivalent:
\begin{enumerate}
\item[(1)] $X$ has strong measure zero.
\item[(2)] $T(X)$ satisfies the selection hypothesis $S_1(\open,\dense)$.
\item[(3)] ONE has no winning strategy $\gone^{\omega}(\open,\dense)$ played on the space $T(X)$.
\end{enumerate}
One can also show that TWO has a winning strategy in $\gone^{\omega}(\open,\dense)$ if, and only if, $X$ is countable.
Thus, by choosing $X\subseteq \reals$ appropriately we obtain a compact ${\sf T}_2$ space which is not weakly Rothberger. 

{\flushleft{\bf F:} Non-Lindel\"of ${\sf T}_2$ space for which TWO has a winning strategy in the game $\gone^{\omega}(\open,\almost)$.}\\
This example demonstrates that none of the implications from the left panel of Figure \ref{cshl} to its middle panel is reversible.\\
The following is Example 3.3 of \cite{PS}. The space $X$ is a subset of $\reals\times\reals$ with a special topology: First, choose a subset $\{x_{\alpha}:\alpha<\omega_1\}$ of $\aleph_1$ distinct elements of the set of nonnegative real numbers. 
\[
  X = \{(x_{\alpha},m): m \mbox{ an integer larger than }-2 \mbox{ and }\alpha<\omega_1\} \cup\{(-1,\,-1)\}.
\]
For convenience we also define:
\[
  A = \{(x_{\alpha},-1):\alpha<\omega_1\}
\] 
and 
\[
  Y = \{(x_\alpha,n): 0\le n<\infty \mbox{ and }\alpha<\omega_1\}.
\]
Topologize $X$ as follows:
Declare each element of $Y$ to be an isolated point; for each $\alpha<\omega_1$ and $n<\omega$ the neighborhood $U_{n,\alpha}$ of $(x_{\alpha},-1)$ is the set $\{(x_{\alpha},-1)\}\cup\{(x_{\alpha},m):m\ge n\}$. Finally, for each $\alpha$ the neighborhood $V_{\alpha}$ of $(-1,1)$ is the set $\{(-1,-1)\}\cup\{(x_{\beta},n):\beta>\alpha,\, -1<n<\omega\}$.

Since the uncountable subset $A$ of $X$ is a closed and discrete subset of $X$, $X$ is not Lindel\"of. Also, for a fixed $\alpha$, for neighborhood $V_{\alpha}$ of $(-1,-1)$ we see that $\overline{V}_{\alpha}$ contains all but countably many elements of $X$. Thus, TWO wins $\gone^{\omega}(\open,\almost)$ as follows: In the first inning TWO chooses the set $T_1$ from the open cover $O_1$ provided by ONE so that $T_1$ contains a neighborhood of $(-1,-1)$ of the form $V_{\alpha}$. Then in the remaining innings TWO makes sure to cover the at most countably many points in the set $X\setminus\overline{T}_1$. 

Moreover, the point $(-1,1)$ does not have a countable neighborhood base, and so this space is not first countable. We now show that this example does meet condition (1) of Theorem \ref{twowinsalmostrothb}: Let $D$ be a countable subset of $X$. We may assume that $(-1,-1)$ is an element of $D$. Fix a $\beta<\omega_1$ such that $D\cap \{(x_{\alpha},n):-1\le n<\omega\} = \emptyset$. Then the neighborhood assignment 
\[
  \{V_{\beta+\omega}\}\bigcup\{\{(x_{\alpha},n)\}: -1<n \mbox{ and } (x_{\alpha},n)\in D\} \bigcup \{U_{0,\alpha}:(x_{\alpha},-1)\in D\}
\]
witnesses the failure of condition (1) of Theorem \ref{twowinsalmostrothb}.

\section{Acknowledgements}

The research for the results reported in this paper occurred while Dr. Pansera was visiting the Department of Mathematics at Boise State University. The Department's hospitality and support during this visit is gratefully acknowledged.



\begin{thebibliography}{88}

\bibitem{LBMS} L. Babinkostova and M. Scheepers, \emph{Combinatorics of open covers (IX): Basis properties}, {\bf Note di Matematica} 22:2 (2003/2004), 167 - 178.

\bibitem{BPSProducts} L. Babinkostova, B.A. Pansera and M. Scheepers, \emph{Weak covering properties and selection principles}, in preparation. 

\bibitem{Bo} E. Borel, \emph{Sur la classification des ensembles de mesure nulle}, \textbf{Bull. Sot. Math. France 47
(1919)}, 97--125.

\bibitem{CHN} W.W. Comfort, N. Hindman and S. Negrepontis, \emph{F-spaces and their products with P-spaces},
\textbf{Pacific Journal of Mathematics} 28 (1969), 489--502.

\bibitem{Daniels} P. Daniels, \emph{Pixley-Roy spaces over subsets of the reals}, {\bf Topology and its Applications} 29 (1988), 93--106.

\bibitem{E} R. Engelking, \emph{General Topology}, 2nd Edition, Sigma Ser. Pure Math., Vol. 6, Heldermann, Berlin,     1989.

\bibitem{F} Z. Frolik, \emph{Generalizations of compact and Lindel\"{o}f spaces}, \textbf{Czechoslovak Math. J. 9 (84) (1959)}, 172--217 (Russian).

\bibitem{Galvinpog} F. Galvin, \emph{Indeterminacy of the point-open game}, {\bf Bull. Acad. Polon. Sci.} 26 (1978), 445 - 449.

\bibitem{Isaac} I. Gorelic, \emph{The Baire category theorem, and forcing large Lindel\"of spaces with points ${\sf G}_{\delta}$}, {\bf Proceedings of the American Mathematical Society} 118 (1993), 603 - 607.

\bibitem{Hurewicz} W. Hurewicz, {\em \"Uber eine Verallgemeinerung des Borelschen Theorems}, {\bf Mathematische Zeitschrift} 24 (1925), 401 -- 421.

\bibitem{Jech} T.J. Jech, \emph{Multiple Forcing}, {\bf Cambridge Tracts in Mathematics} 88, 1987.

\bibitem{COC2} W. Just, A. W. Miller, M. Scheepers, P. J. Szeptycki, \emph{The combinatorics of open covers (II)}, \textbf{Topology and its Applications} 73 (1996), 241--266

\bibitem{Kada} M. Kada, \emph{Preserving the Lindel\"{o}f property under forcing extensions}, {\bf Topology Proceedings} 38 (2011),  237--251.

\bibitem{K} D. Kocev, \emph{Almost Menger and related spaces}, \textbf{Matemati\v cki Vesnik}
61 (2009), 173--180.

\bibitem{Ko2} Lj.D.R. Ko\v{c}inac, \emph{Star-Menger and related spaces} II, \textbf{Filomat} 13 (1999), 129--140.

\bibitem{P} B.A. Pansera, \emph{Weaker forms of the Menger property}, \textbf{Quaestiones Mathematicae} 34 (2011), 1--9.

\bibitem{pawlikowski} J.Pawlikowski, \emph{Undetermined sets of point-open games}, {\bf Fundamenta Mathematicae} 144 (1994), 279 -- 285.

\bibitem{PW} J.R. Porter, R.G. Woods, \emph{Extensions and absolutes of Hausdorff spaces}, \textbf{Springer-Verlag}, 1988.

\bibitem{R} F. Rothberger, \emph{Eine Versch\"{a}rfung der Eigenschalft \textsf{C}}, \textbf{Fundamenta Mathematicae 3 (1938)}, 50--55.

\bibitem{COC1} M. Scheepers, \emph{Combinatorics of open covers I: Ramsey theory}, \textbf{Topology and its Applications 69(1996)}, 31--62.

\bibitem{COC4} M. Scheepers, \emph{Combinatorics of open covers (IV): subspaces of the Alexandroff double unit interval}, \textbf{Topology and its Applications 83(1998)}, 63--75.

\bibitem{COC5} M. Scheepers, \emph{Combinatorics of open covers (V): Pixley-Roy spaces of sets of reals, and $\omega$-covers}, \textbf{Topology and its Applications 102(2000)}, 13--31.

\bibitem{MSLusin} M. Scheepers, \emph{Lusin sets}, \textbf{Proceedings of the American Mathematical Society 197(1999)}, 251--257.

\bibitem{MST} M. Scheepers, \emph{A direct proof of a theorem of Telg\'arsky}, {\bf Proceedings of the American Mathematical Society} 123:11 (1995), 3483 - 3485.

\bibitem{MStightness} M. Scheepers, \emph{Remarks about countable tightness}, ArXiv 1201.4909

\bibitem{MSRothbgroup} M. Scheepers, \emph{Rothberger bounded groups and Ramsey theory}, {\bf Topology and its Applications} 158 (2011), 1575 - 1583.

\bibitem{ST} M. Scheepers and F.D. Tall, \emph{Lindel\"of indestructibility, topological games and selection principles}, {\bf Fundamenta Mathematicae} 210 (2010), 1 - 46

\bibitem{SS} L.A. Steen, J.A. Seebach, \emph{Counterexamples in Topology}, \textbf{Dover Publications, New York}, 1995.

\bibitem{PS} P. Stoyanova, \emph{A comparison of Lindel\"of type covering properties of topological spaces}, {\bf Rose-Hulman Undergraduate Mathematics Journal} 12:2 (2011), 163 - 204.

\bibitem{Tall} F.D. Tall, \emph{On the cardinality of Lindel\"of spaces with points ${\sf G}_{\delta}$}, {\bf Topology and its Applications} 63 (1995), 21 - 38.

\bibitem{Telgarsky} R. Telg\'arsky, \emph{On games of Topsoe}, {\bf Math. Scand.} 54 (1984), 170 - 176.

\bibitem{Tkachuk} V.V. Tkachuk, \emph{Some new versions of an old game}, {\bf Commentationes Mathematicae Universitatis Carolinae} 36:1 (1995), 177 - 196.

\bibitem{WD} S. Willard, U.N.B. Dissanayake,\emph{The almost Lindel\"{o}f degree}, \textbf{Canad. Math. Bull., 27
(4) (1984)}, 452–-455 .

\end{thebibliography}
\end{document}